\documentclass[aps,pre,twocolumn]{revtex4}
\usepackage{amsmath}
\usepackage{amsfonts}
\usepackage{amssymb}
\usepackage{graphics}
\usepackage[pdftex]{graphicx}

\def\bee{\begin{equation}}
\def\eee{\end{equation}}
\def\Li{{\rm Li}}

\def\DHLhksqrt#1#2{\setbox0=\hbox{$#1\sqrt{#2\,}$}\dimen0=\ht0
\advance\dimen0-0.2\ht0
\setbox2=\hbox{\vrule height\ht0 depth -\dimen0}%
{\box0\lower0.4pt\box2}}
\begin{document}

\title{ Nearest neighbor spacing distribution  of prime numbers\\ and quantum chaos}

\author{Marek Wolf}
\affiliation{Cardinal  Stefan  Wyszynski  University, Faculty  of Mathematics and Natural Sciences. College of Sciences\\
ul. W{\'o}ycickiego 1/3, PL-01-938 Warsaw,   Poland, e-mail:  m.wolf@uksw.edu.pl}

\begin{abstract}
We give heuristic arguments and  computer results to support the hypothesis that, after
appropriate rescaling, the statistics of spacings between adjacent prime numbers follows
the Poisson distribution. The  scaling  transformation removes the oscillations in the NNSD of primes.  These oscillations
have  the very  profound period  of  length six.  We also calculate the spectral rigidity $\Delta_3$ for prime numbers
by two methods. After suitable averaging one of these methods gives the Poisson dependence $\Delta_3(L)=L/15$.
\end{abstract}

\maketitle

\bigskip\bigskip\bigskip

%Key words: {\it Prime numbers, gaps between primes, Hardy and Littlewood conjecture }\\  PACS numbers: {11A99, 11A41}

\bibliographystyle{apsrev} %amsplain} % unsrt }

\section{Introduction}

The primes numbers  often  provided  a toy model for some physical ideas in the past. For example
in \cite{Wolf-1989} the multifractal formalism was applied to prime numbers, in \cite{Gamba-Hernando}
the appropriately defined Lyapunov exponents for the distribution of primes were calculated
numerically.  In the  paper  \cite{Wolf-1997}  it was  shown that the distribution
of prime numbers displays the $1/f$ noise, while in \cite{Lan-Yong-2004} the noise $1/f^2$ was found in the difference
between the prime-number counting  $\pi(x)$  function and Riemann's function $R(x)$. In \cite{Billingsley1973} and  \cite{Wolf1998}
random walks on primes numbers were defined. In \cite{Bonanno-2004} an attempt to construct the dynamical model for
prime numbers was taken  and computable information content as well as entropy information of the set of  prime numbers were
calculated.

The prime numbers can be regarded as eigenvalues of some quantum hamiltonian.   The problem of construction of a simple
one--dimensional Hamiltonian whose spectrum coincides with the set of primes was considered in \cite{Mussardo},
\cite{Sekatskii-2007}, \cite{Schumayer-et-al},   see  also  review  \cite{Rosu-2003}.
Then it is natural to investigate the spacings between prime numbers, i.e. in physical language
the nearest neighbor spacing distribution (NNSD).  Several authors have undertaken a study of this problem
in the past, see \cite{Liboff_Wong-1998},  \cite{Timberlake-2006}, \cite{Timberlake-2007}.  Below we will treat prime
numbers as the energy levels and we will apply methods used to describe statistical properties of discrete  spectra.
Let the quantum  system possess  the discrete spectrum $E_1, E_2, \ldots$ and
let $N(E)=\sum_{n}\Theta(E-E_n)$ ($\Theta$ is a unit step function) denote the function counting the number of energy levels smaller than $E$.
Usually spectral staircase $N(E)$ can be split into the ``smooth'' $\overline{N}(E)$ and fluctuating (oscillating) $\widetilde{N}(E)$
parts.  For example,  for a large class of differential operators on $d$  dimensional bounded manifold $\Omega \subset \mathbb{R}^d$
the Weyl's law
\bee
\overline{N}(E) \sim \frac{\mbox{vol}(\Omega)}{(2\pi)^d}E^{d/2},
\eee
holds,  see e.g.  \cite[Ch.1]{Arendt-2009} .

Given the spectrum $E_1, E_2, \ldots$ the statistics of normalized and dimensionless (``unfolded'' spectrum, see e.g.
\cite[Sect.4.7]{Haake}) gaps between two consecutive energy levels
$s_n =(E_{n+1}-E_n))/\overline{d}(E)$, where $\overline{d}(E)$ is the mean distance
between energy levels up to $E$,  was extensively studied in the past. For general systems $E_{n+1}-E_n$ are arbitrary real
numbers and histogram of the level spacings  $s_n$  is built.  It is well known, that level--spacing
distributions of quantum systems can be grouped into a few universality  classes connected with the symmetry
properties of the hamiltonians: Poisson distribution (i.e. $e^{-s}$)  for systems with underlying regular classical dynamics,
Gaussian orthogonal ensemble (GOE, also called the Wigner--Dyson distribution) --- hamiltonians invariant under time reversal,
Gaussian  unitary ensemble (GUE) --- not invariant under time reversal  and Gaussian symplectic ensemble (GSE)
for half-spin systems with time reversal symmetry. There are many reviews on these topics, we cite here
\cite{Mehta}, \cite{Haake}, \cite{Weidenmuller-Mitchell}.   % {Mehta, Haake, Weidenmuller-Mitchell}.% ,

%  par. 3.2.4  (od str. 110)  w  Stockmann, Quantum Chaos_An_Introduction.pdf

There is some confusion  regarding the proper statistics of the gaps between consecutive
primes: in \cite{Liboff_Wong-1998} it was claimed that NNSD of primes follows GOE distribution, while
in  \cite{Timberlake-2006, Timberlake-2007},  the  possibilities of GOE, Poisson and exotic
Berry-Robnik  \cite{Berry-Robnik} distribution were investigated.  Liboff  and  Wong   have obtained Wigner
distribution and level repulsion
for NNSD of primes  by artificially including the gaps 0 (no degeneracy --- all primes are different) and 1, see
\cite[p.3113]{Liboff_Wong-1998}. The gap 1 appears only once
between 2 and 3 and should be skipped in the wake of  infinity of primes. There is a very often reproduced  % widespread
figure showing some  typical spectra (see \cite[Fig. 1.2]{Mehta}, \cite[Fig.3]{Weidenmuller-Mitchell}, \cite[Fig. I.8]{Bohigas-1984},
\cite[front figure]{Cipra-1999}, \cite[p. 32]{Terras-2011}): random levels with no correlations (Poisson series),
sequence of prime numbers, resonance levels of erbium 166 nucleus,  the  energies a free particle in the Sinai billiard,
nontrivial zeros of the Riemann zeta function. In \cite[p. 10]{Mehta}  it is stated that ``case of prime numbers $\ldots$
are far from either regularly spaced uniform series or the completely random Poisson  series with no correlations''.

It is the purpose of this paper ``to settle  once and for ever'' that NSDD of primes follows the Poisson distribution.
The next Section II is devoted to  this problem. In \cite{Berry-1985}  M.V.  Berry has  calculated spectral rigidity
$\Delta_3$ for zeros of the Riemann  zeta function  and in Sect.III we will study spectral rigidity for prime numbers.  % We will exploit  two methods of calculating   $\Delta_3$

\begin{figure*}[th]
\includegraphics[width=\textwidth, angle=0]{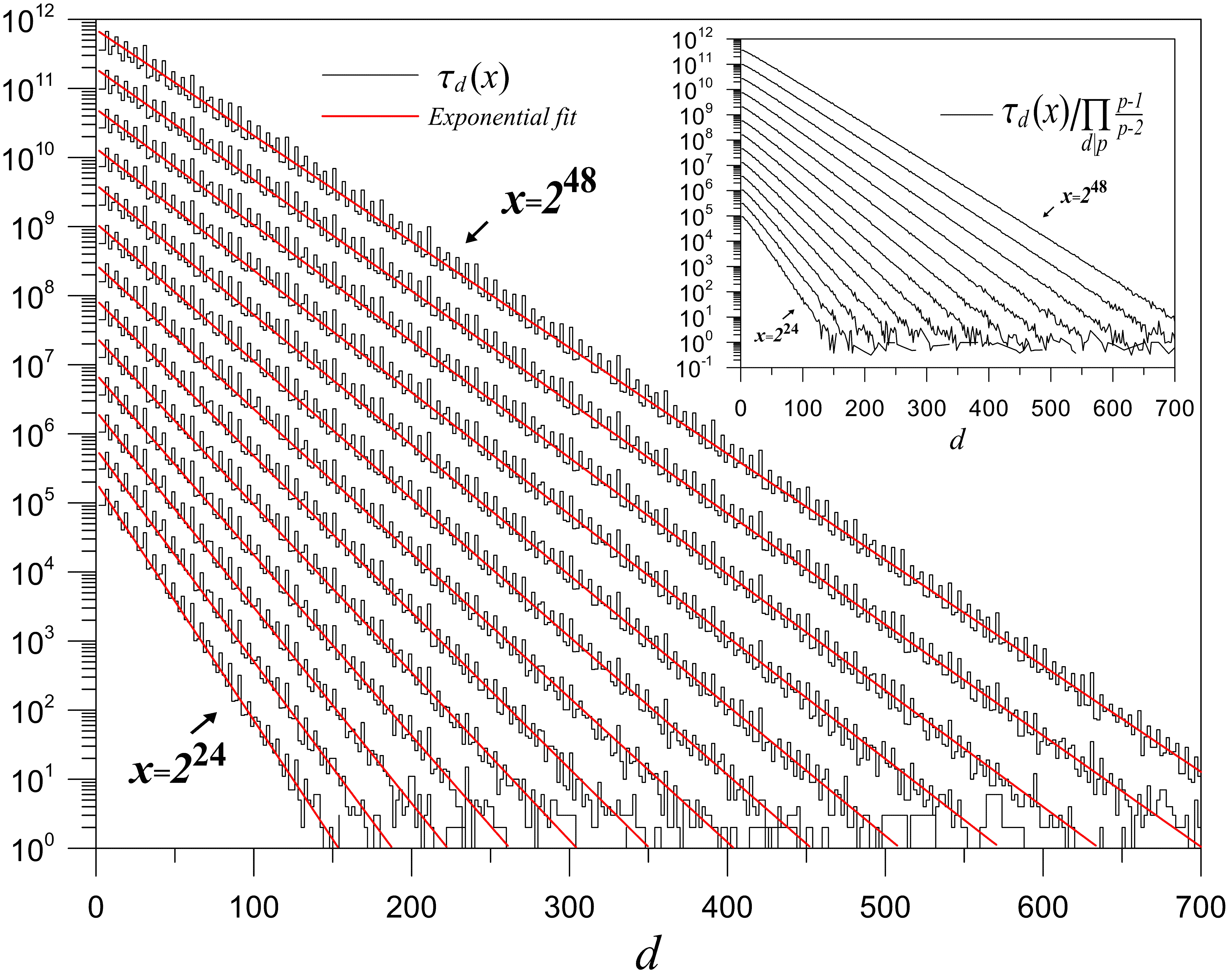}
\caption{Plots of $\tau_d(x)$ for $x=2^{24}, 2^{26}, \ldots, 2^{46}, 2^{48}$.
The histogram step widths are $2$; because $\tau_2(x)\approx \tau_4(x)$, therefore the visible step for $d=2, 4$ has width $4$.
In red  exponential fits $a(x)e^{-db(x)}$  are plotted. In the inset the plots of $\tau_d(x)/P(d)$ are shown.}
\label{fig-hist}
\end{figure*}

%!!!!!!!!!!!!!!!   Aurich  Steiner Staircase functions, spectral rigidity, and a rule for quantizing chaos PRA 45 p583_1.pdf   *********

\section{NNSD for prime numbers}

In the case of primes numbers all gaps $d_n=p_{n+1}-p_n$ (except the first pair of   primes $p_1=2, p_2=3$) are even
integers $2,4, 6, \ldots$. These spacings are dimensionless and we will not perform unfolding for time being (see next Section) ---
the usual (\ref{unfolding-PNT}) unfolding obscures analysis of the oscillations present in the NNSD between  original primes.
Let $\tau_d(x)$ denote a number of pairs
of {\it consecutive} primes smaller than a given bound $x$ and separated by $d$:
\bee
%\tau_d(x)= \{{\rm number~ of~ pairs~~} p_n,  p_{n+1} < x,~~{\rm  with}~ d=p_{n+1}-p_n\}.
%\tau_d(x)= \{{\rm \sharp~ pairs~~} p_n,  p_{n+1} < x,~~{\rm  with}~ d=p_{n+1}-p_n\}.
\tau_d(x)= \sharp\{ p_n,  p_{n+1} < x,~~{\rm  with}~ p_{n+1}-p_n=d\}.
\label{definition}
\eee
For odd $d=2k+1$ we  supplement this definition by putting $\tau_{2k+1}(x)=0$.

In 1922   G. H. Hardy and J.E. Littlewood in the famous paper  \cite{Hardy_and_Littlewood}
have proposed 15 conjectures. The  conjecture B of their paper states that
there are infinitely many primes pairs $(p, p^\prime)$, where $p^\prime = p+d$,
for every even $d$. If $\pi_d(x)$ denotes the number of prime pairs differing
by $d$ and  less than $x$, then
\bee
\pi_d(x)\sim  C_2  \prod_{p\mid d} \frac{p-1}{p-2}~ \frac{x}{\ln^2(x)}.
\label{H-L}
\eee
Here $C_2 \equiv 2\prod_{p > 2} \biggl( 1 - \frac{1}{(p - 1)^2}\biggr) =1.32032\ldots $ is
called the ``twins constant''.

In the middle of 2013 the major step towards the proof of the conjecture B was made: Yitang Zhang
has submitted to Annals of Mathematics the paper in which he proved unconditionally that
$\liminf_{n\to \infty}\, (p_{n+1} - p_n) < 7 \times 10^7$, see e.g. \cite{Zhang-Nature}.  Very soon this bound was lowered many
times by  mathematicians  and present record is $\liminf_{n\to \infty}\, (p_{n+1}-p_n)\le 600$ and was obtained
by J. Maynard \cite{Maynard-2013}.

The conjecture B of G. H. Hardy and J.E. Littlewood  gives the number of pairs of primes
not necessarily successive  and we would like to  stress that in  (\ref{definition})  $\tau_d(x)$  denotes
number of pairs of {\it consecutive} primes $p_n, p_{n+1}$ with difference $p_{n+1} - p_n = d$.
The pairs of primes separated by $d=2$ and $d=4$ are special as they
always have to be consecutive primes (with the exception of the pair
(3,7) containing 5 in the middle): in the triple of integers $2k+1, 2k+3, 2k+5$ the middle
$2k+3$ has to be divisible by 3 if $2k+1, 2k+5$ are prime (in particular not divisible by 3). For $d=6$ (and larger $d$)
we have $\pi_6(x)>\tau_6(x)$, for example  $(5, 7, 11), (7, 11, 13), (11, 13, 17), \ldots$.
From the conjecture B of G. H. Hardy and J.E. Littlewood  \cite{Hardy_and_Littlewood}
it follows that the number of gaps $d=2$ (``twins'') is approximately equal
to the number of gaps $d=4$ (``cousins''): $\pi_2(x)\equiv \tau_2(x)\approx \pi_4(x)\equiv \tau_4(x)$, see also \cite{Wolf1998}.
For $d\geq 6$ in  \cite{Wolf-1999}  we have conjectured that

\begin{widetext}
\bee
\tau_d(x)\sim  C_2 \frac{\pi^2(x)}{x} \prod_{p \mid d, p > 2}\frac{p - 1}{p - 2} e^{- {d \pi(x)/x}}~~ {\rm for} ~ d\geq 6,~~~~\tau_2(x)\big(\approx \tau_4(x)\big)\sim
C_2 \frac{\pi^2(x)}{x} \approx C_2 \frac{x}{\ln^2(x)}.
\label{main}
\eee
\end{widetext}

\noindent Here  $\pi(x)=\sum_{n} \Theta(x-p_n)$ denotes the number of primes up to $x$  and by
the Prime Number Theorem (PNT) is very well approximated by the logarithmic integral
\[
\pi(x)\sim  {\rm Li}(x) \equiv \int_2^x \frac{du}{\ln(u)}.
\]
Integration by parts gives the asymptotic expansion which should be cut at the term $n_0=\lfloor \ln(x)\rfloor$:
\bee
{\rm Li}(x) = % \frac{x}{\ln x} \sum_{k=0}^{\infty} \frac{k!}{(\ln x)^k} =
\frac{x}{\ln(x)}  + \frac{x}{\ln^2(x)} + \frac{2!x}{\ln^3(x)} + \frac{3!x}{\ln^4(x)} +\cdots .
\label{PNT}
\eee
There is a series giving $\Li(x)$ for  all $x>2$ and  quickly convergent which has $n!$ in denominator and $\ln^n(x)$
in nominator instead of opposite order in (\ref{PNT})
(see \cite[Sect. 5.1]{abramowitz+stegun})
\begin{equation}
\Li (x) =  \gamma +
\ln \ln(x) + \sum_{n=1}^{\infty} {\ln^{n}(x)\over n \cdot n!}
\quad {\rm for} ~ x > 1 ~ ,
\label{Li-series}
\end{equation}
Here $\gamma=0.577216...$ is the Euler-Mascheroni constant.

Putting in (\ref{main})  $\pi(x)\sim  {x}/{\ln(x)}$  the compact formula
expressing $\tau_d(x)$ by  explicitly known functions
\bee
\tau_d(x)\sim  C_2 \frac{x}{\ln^2(x)} \prod_{p \mid d, p > 2}\frac{p - 1}{p - 2} e^{- {d /\ln(x)}}~~
\label{main2}
\eee
is obtained.  Comparing it with the original Hardy--Littlewood conjecture (\ref{H-L})  we obtain that
the number $\tau_d(x)$ of {\it successive} primes $(p_{n+1}, ~p_n)$ smaller
than $x$ and of the difference $d~(=p_{n+1}-p_n)$ is diminished by the factor $\exp(-d/\ln(x))$
in comparison with the number of {\it all} pairs of primes $(p, p^\prime)$ apart in
the distance $d=p^\prime-p$:
\bee
\tau_d(x) \sim  \pi_d(x) e^{-d/\ln(x)}  ~~~~~{\rm for~} d\geq 6.
\label{relation}
\eee
The expression (\ref{main2})  for $\tau_d(x)$ was proved (in slightly different form required by the precision of
the formulation of the theorem) under the assumption of the conjecture B
of Hardy--Littlewood by  D. A. Goldston and A. H. Ledoan \cite{Goldston-Leodan-2012} in 2012.

\begin{figure*}[ht]
\begin{minipage}[ht]{0.5\textwidth}
\includegraphics[width=\textwidth, angle=0]{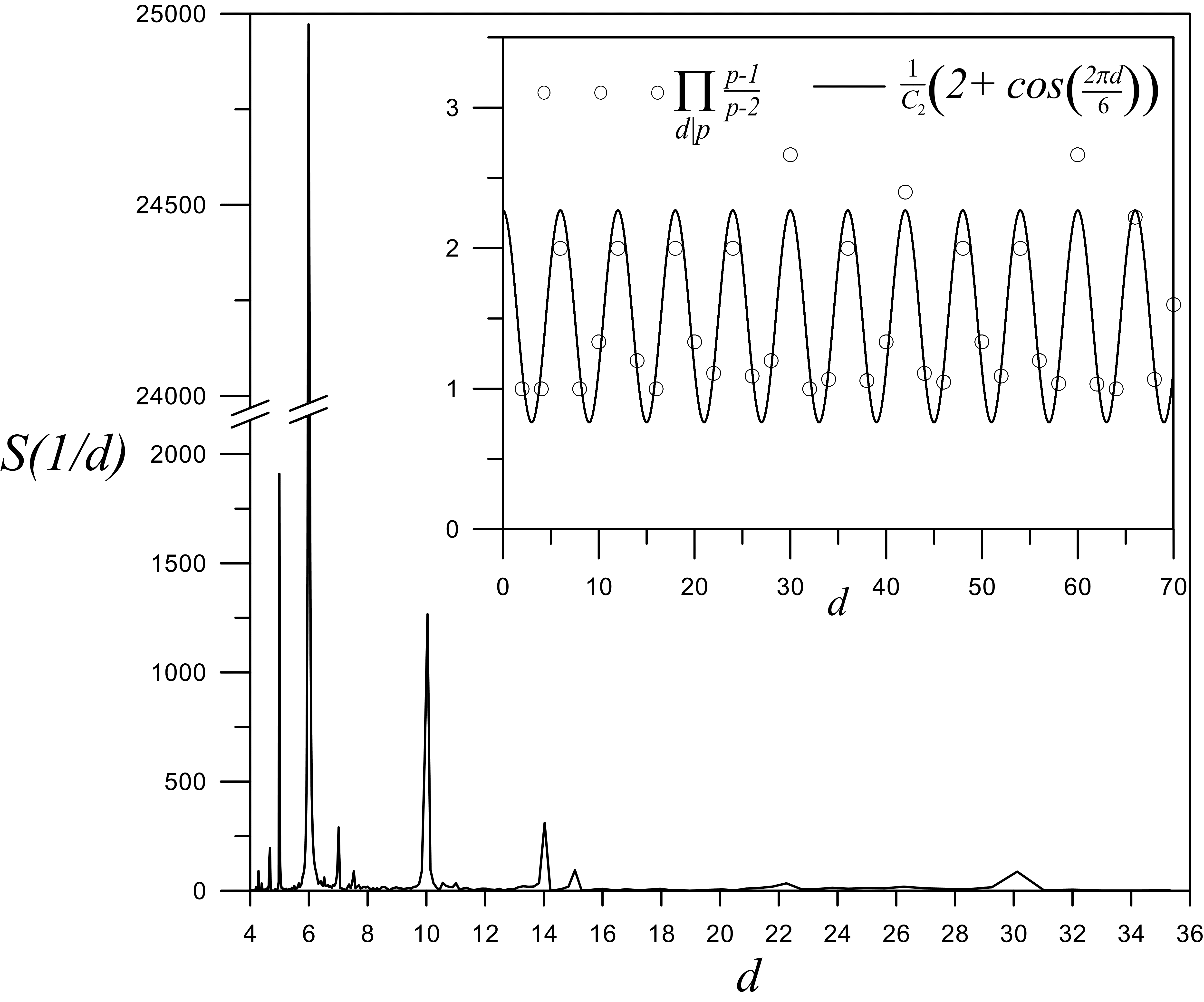} \\%{NNSD_a.png} \\
\caption{The plot of power spectrum $S(f)$  calculated from $M=2^{10}=1024$  values of $P(d)$ plotted versus
$1/f$ to show main periods of $P(d)$. The $y$ axis  was broken  to make visible peaks at $d\neq 6$.
 In the inset the plots of $P(d)$ and the approximation (\ref{ansatz2}) are presented. }
\label{fft}
\end{minipage}
\hspace{0.03\textwidth}
\begin{minipage}[ht]{0.45\textwidth}
\centering
\includegraphics[width=\textwidth,angle=0]{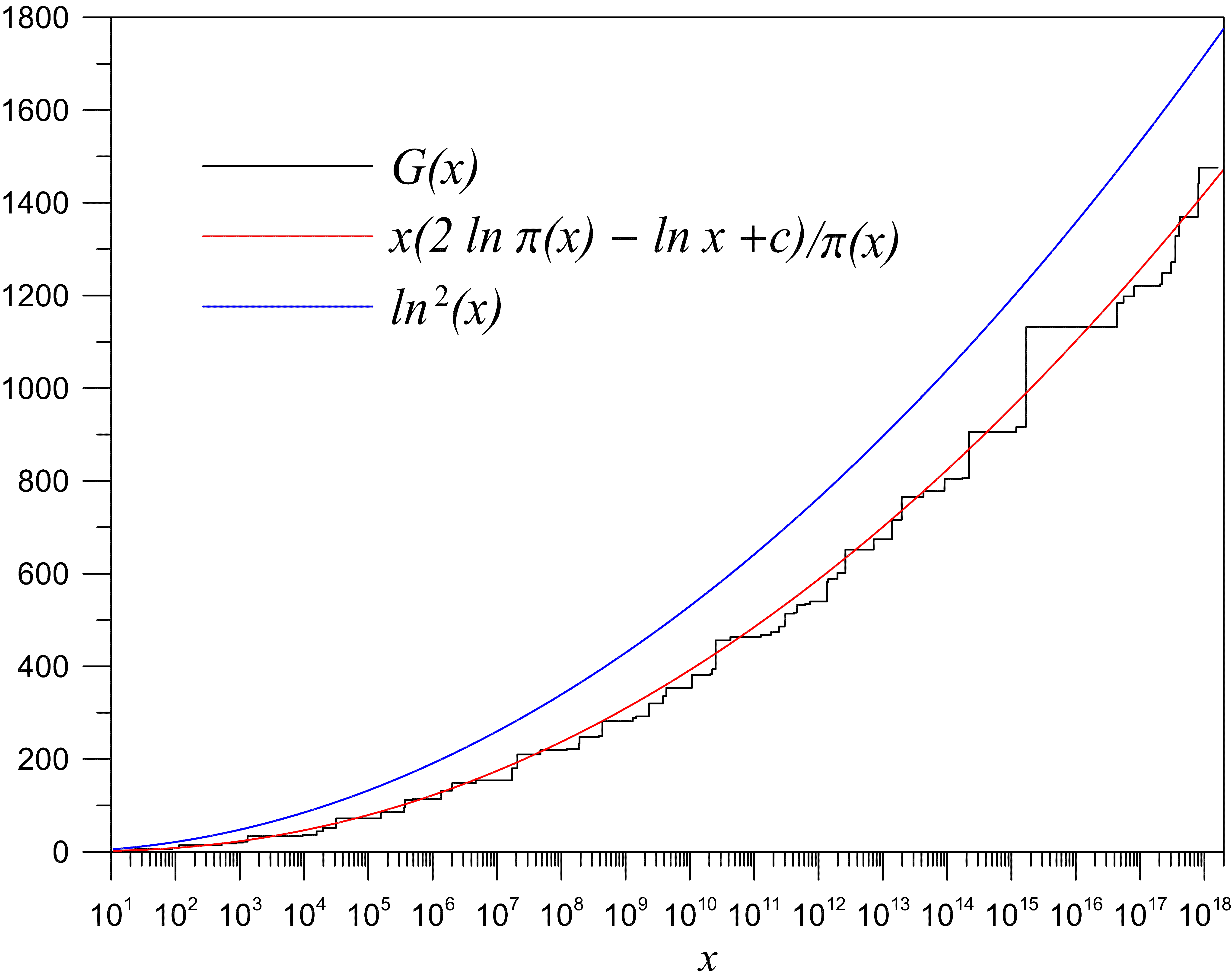}
\caption{The comparison of $G(x)$ obtained from the computer search (up to $x=2^{48}$ we have used our own data, for
larger $x$ we took data from the web pages \cite{web-primes}). For the plot of (\ref{d_max}) we have used the tabulated
values of $\pi(x)$ available at \cite{web-primes}. The plot of the Cramer conjecture is also presented. }
\label{fig-G}
\end{minipage}
\end{figure*}

During over a seven months long run of the computer
program we have collected  the values of $\tau_d(x)$ up to $x=2^{48}\approx 2.8147\times 10^{14}$.
The data representing the
function $\tau_d(x)$ were stored at values of $x$ forming the geometrical
progression with the ratio 2  at $x=2^{15}, 2^{16}, \ldots, 2^{47},
2^{48}$. Such a choice of the intermediate thresholds as powers of 2
was determined by the employed computer program in which  the primes were coded
as bits. The data is available for downloading from  http://pracownicy.uksw.edu.pl/mwolf/gapstau.zip.
The resulting curves are plotted in Fig.\ref{fig-hist}.  Characteristic oscillating pattern of points is
caused by the product
\bee
P(d)\equiv  \prod_{p \mid d, p > 2} \frac{p - 1}{p - 2}
\label{product}
\eee
appearing in (\ref{main}), see inset in Fig. \ref{fig-hist}. This product for the first time appeared in the paper of Hardy
and Littlewood \cite{Hardy_and_Littlewood} and it has local maxima for $d$ equal to  the products of consecutive primes (``primorials'',
i.e. factorials over primes $2\cdot 3 \cdot 5 \ldots \cdot p_n\equiv p_n\sharp$).
Clearly visible  in Fig. \ref{fig-hist} are oscillations of the period
$6=2\times 3$ with overimposed  higher harmonics $30=2\times 3\times 5$  and  $210=2\times 3\times 5\times 7$, i.e. when $P(d)$
has local maxima  $P(6)=2, ~P(30)=8/3=2.666\ldots ~ P(210)=16/5=3.2$  (local minima are 1 and they correspond to $d=2^m$).
We have  performed the discrete  Fourier Transform  of $P(d)$, i.e. we calculated numerically
\bee
\widetilde{P}\left(\frac{n}{2M}\right)= \sum_{k=0}^{M-1} P(2k) e^{2\pi kn/M},
\eee
where $n=0, 1, 2, \ldots, M-1$ and $n/2M$ plays the role of discrete frequency.
Having $\widetilde{P}(f)$ we can calculate the power spectrum density  $S(f) = \mid\!\! \widetilde{P}(\frac{n}{2M})\!\! \mid^2$.
The  large value of $S(f)$ at some frequency  $f$ means that the dependence
of $P(d)$ on $d$  has some harmonic component of the period  $T=1/f$. Thus in the Fig.\ref{fft} we have
plotted $S(f)$  versus $1/f=d$ to show main periods  $5, 6=2\times 3, 10=2\times 5, 14=2\times 7, 30=2\times 3\times 5\ldots $
of $P(d)$.  These oscillations are the reason why the  Poisson distribution was
not attributed to NNSD of primes in the past:  e.g. $P(2)=P(4)=1$  while $P(6)=2$ and  the plot should be
made with logarithmic scale on the $y$ axis  to suppress these oscillations.

In \cite{Bombieri} E. Bombieri and H. Davenport have proved that:
\bee
\sum_{k=1}^n \prod_{p\mid k,p>2}\frac{p-1}{ p-2} =
\frac{n}{\prod_{p > 2}( 1 - \frac{1}{(p - 1)^2})} + \mathcal{O}(\ln^2(n));
\label{srednia}
\eee
i.e. in the limit $n\to \infty$ the number  $2/C_2$
is the arithmetical average of the product $\prod_{p\mid k} \frac{p-1}{p-2}$.
The main period of oscillations is $6$  hence  we can write:
\bee
P(d)= \prod_{p \mid d, p > 2} \frac{p - 1}{p - 2} \approx \alpha + \beta\cos\left(\frac{2\pi d}{6}\right).
\label{ansatz}
\eee
The numerical value of $\alpha$ is  equal to $2/C_2$ to reproduce  the average value of $P(d)$ in (\ref{srednia}).
It  can be explained by taking into account that
$\cos(2\pi2k/6)=1$   while  $\cos(2\pi(2k+2)/6)=\cos(2\pi(2k+4)/6)=-\frac{1}{2}$  and hence by the equation
\bee
\lim_{n\rightarrow \infty}\frac{1}{n} \sum_{k=1}^n \Big( \alpha + \beta\cos \bigg(\frac{2\pi 2k}{6}\bigg)\Big) = \alpha
\eee
the value of the parameter $\beta$ does not contribute to the average of r.h.s.  of (\ref{ansatz}).
Thus from (\ref{srednia})  we have $\alpha=2/C_2\approx 1.5147801281$.  Requiring, that the combination
$\alpha+\beta\cos(2\pi d/6)$ for $d=6$  takes the value  2 times   larger than  for  $d=2$ and $d=4$:  $\alpha+\beta=2(\alpha-\beta/2)$
gives $\beta=\alpha/2\approx 0.75739$. Fitting of the parameters $\alpha$ and $\beta$ can be done also numerically by
standard General Linear Least Squares, see e.g.\cite{recipes}. We  have used the procedure
{\it lfit}  from \cite{recipes}   with 2500000 numbers of points:  for $d=2, 4, 6, \ldots, 5000000$. The output
of the computer run was:  $\alpha= 1.51478\approx 2/C_2,~~~ \beta= 0.75471\approx 1/C_2$. Hence we propose the compact formula
(see inset in Fig. \ref{fft}):
\bee
P(d)= \prod_{p \mid d, p > 2} \frac{p - 1}{p - 2} \approx \frac{1}{C_2}\left(2+\cos\left(\frac{2\pi d}{6}\right)\right),
\label{ansatz2}
\eee
which allows to substitute  for  $P(d)$ an expression more amenable for
algebraic manipulations.   Such an approximation may be relevant for calculations
of correlations functions for zeros of the Riemann zeta  function, where sums involving product $P(d)$ appear
very often \cite{Keating-1993}. It turns out, that $\cos(2\pi d/6)$ takes for even $d$ only two values:
$-1/2$ for $d=6k+2$ and  $6k+4$, and 1 for $d=6k$. Because $d$ and $d^2$ have the same prime divisors it follows
that $P(d^2)=P(d)$.   The same relation is also obeyed by the  approximation (\ref{ansatz2}) because $(6k+2)^2=6k^\prime+4$ and
$(6k+4)^2=6k^{\prime\prime}+4$ and the square of the $d=6k$ is obviously again a number of the same form.

The smallest gap between adjacent primes is 2 (twin primes), while the maximal gap $G(x)={\max_{p_n<x}}~(p_n-p_{n-1})$
grows with $x$. We can obtain the formula for $G(x)$ from (\ref{main}) assuming that the largest gap  up to $x$
between two consecutive  ``levels'' $p_{n+1}-p_n$  appears only once: $\tau_{G(x)}(x)=1$. Skipping the
oscillating term $P(d)$, which is very often close to 1,
we get for $G(x)$ the following estimation expressed directly by $\pi(x)$:
\bee
G(x) \sim \frac{x}{\pi(x)} \big(2\ln( \pi(x)) -\ln(x) + c\big),
\label{d_max}
\eee
where $c=\ln(C_2)= 0.2778769\ldots$.  Substituting here the PNT in the form $\pi(x)\sim x/\ln(x)$ gives the
Cramer's conjecture  \cite{Cramer}  $G(x)\sim \ln^2(x)$ in the limit of large $x$. The maximal gaps $G(x)$ are scattered
chaotically,   the largest currently known gap of 1476 follows the prime 1425172824437699411,  see  \cite{web-primes}.
The comparison of the above formula with real data  is presented in Fig. \ref{fig-G}.

We finish this section recalling the result of P. Gallagher \cite{Gallagher1976}.
He proved,  assuming the special generalization of the $n$-tuple
conjecture of Hardy--Littlewood  (\ref{H-L}), that the fraction of intervals which
contain exactly $k$ primes follows a Poisson distribution. More precisely he proved, that
the number $P_k(h, N)$ of such $n<N$ that the interval $(n , n+h]$   contains exactly $k$ primes
is asymptotically for $N \to \infty$ given by
\[
P_k(h, N)\sim N \frac{\lambda^k e^{-\lambda}}{k!},
\label{Gallagher}
\]
where $\lambda \sim h/\ln(N)$ is a parameter of the Poisson distribution.  In \cite{Kowalski2011} E. Kowalski has generalized
the Gallagher theorem to other families of primes. In particular the numbers of twins, primes of the form $m^2+1$
or Sophie Germain primes (i.e. primes p with 2p + 1  also prime)  in  short intervals are asymptotically Poisson distributed.

\section{Unfolded primes}

\begin{figure}
\includegraphics[width=0.45\textwidth, angle=0]{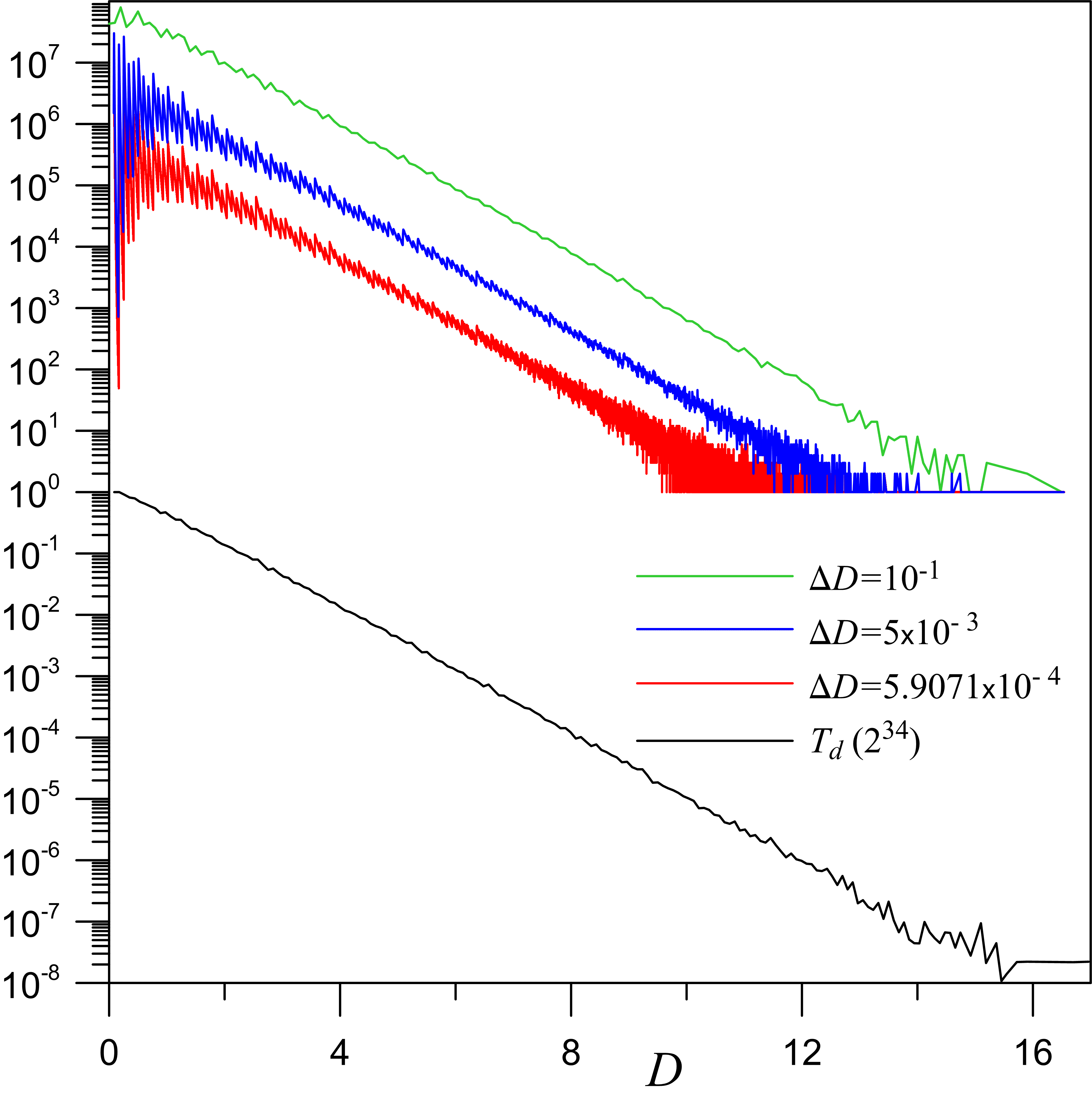}
\caption{The plot of histograms of unfolded spacings $D_n=r_{n+1}-r_n$ where $r_n=\Li(p_n)$ for primes up to $2^{34}=1.72\ldots \times 10^{10}$.
Three widths of bins are used: $\Delta D=0.1$,  $\Delta D=0.001$ and $\Delta D=16.54/28000$.
In black is shown the plot for the unfolding defined by eq. (\ref{eq-scaling}).}
\label{fig-bins}
%\vspace{0.6cm}
\end{figure}

For energy spectrum $E_1, E_2, \ldots$   one  usually performs unfolding to focus on fluctuations around the smooth part of
staircase  and simultaneously to pass to the dimensionless  variables $e_1, e_2, \ldots$   via the definition:
\bee
e_n=\overline{N}(E_n).
\label{unfolding}
\eee
Then the  average spacing between two consecutive $e_n, e_{n+1}$ is  equal to 1 and this procedure removes the
individual properties of a system.  Although primes are dimensionless  we can perform the unfolding using
the definition
\bee
r_n=\Li(p_n).
\label{unfolding-PNT}
\eee
Then the unfolded spacings are $D_n=r_{n+1}-r_n$, writing $p_{n+1}=p_n+d_n$ ($d_n$ are ``pure'' spacings, not unfolded)
and using $\Li(x)\sim x/\ln(x)$ we obtain
\bee
D_n\approx \frac{d_n}{\ln(p_n)+d_n/p_n}
\label{D_n-formula}
\eee
and for large $p_n$ it goes into $D_n = d_n/\ln(p_n)$.  In other words we can say, that the unfolded gaps
(level spacings) between very large  consecutive primes are $D_n=(p_{n+1}-p_n)/\ln(p_n)$.
Because  the average distance between  primes  $(p_{n-1}, p_{n})$  is $\ln(p_n)$ we have  from (\ref{D_n-formula}) for
large $p_n$  that the average spacing between two consecutive $(r_n, r_{n+1})$ is equal to 1,  as it should be for unfolded variables.
The values of $D_n$ are arbitrary real numbers, while $d_n$ assume only even values. For example, for twin primes $p_{n+1}=p_n+2$
the  gap  $d=2$  will be mapped into $D_n\approx 2/\ln(p_n)$  with explicit dependence on $p_n$ and it goes to zero with increasing
$p_n$ (if there are infinity of twins, as it is widely believed).  On the other side the maximal value of $D$ will correspond
to maximal gaps: from (\ref{d_max}) we have that roughly $G(p_n)=\ln^2(p_n)$ and thus the interval of values of $D$ will
span up to approximately $\ln(p_n)$:  the values $d=2,\, 4,\, 6\, \ldots,\, G(x)$ will be mapped onto the interval
$[2/\ln(x), \ln(x)]$. To make the histogram of unfolded spacings  $D_n$ the (arbitrary) size of bin should be chosen.
In this approach the oscillations seen in Fig. \ref{fig-hist} are ``smeared out'' between  different bins  and there
is no possibility to extract them easily from the histogram  of unfolded  gaps $D_n$ --- the behavior  caused by
the product $P(d)$  is  obscured after the change of variables  $d_n\rightarrow  D_n$, see  oscillations with large amplitude
on the red and blue plots in Fig. \ref{fig-bins} --- $D_n$  depends explicitly on the value of $p_{n}$ and is a continuous variable.
In other words  the same bin will contain contribution from different $d_n$ and different
$p_n$  giving the same value of $D_n$ and there is no possibility to
untangle  for unfolded quantities the influence of the oscillations caused by the product (\ref{product}).
We present the results of this procedure for all primes up to $2^{34}=1.718\ldots\times 10^{10}$ in Fig.\ref{fig-bins} for three
choices of the bin size.  The popular choice, used   e.g. by spreadsheet  Excel, is to set the number of bins equal
to the square root of the number of values of binned  variable. In our case $\pi(2^{34})=762939111$,  thus the number of bins
should be approximately 28000. Because the maximal gap up to $2^{34}$ is  $G(2^{34})=382$  and it appears at $p_{486570087}=10726904659$
we get that the maximal value of $D$ is $382/\ln(10726904659)=16.54\ldots$ and the size of bin should be $16.54/28000\approx 0.00059$.
In Fig. \ref{fig-bins}  red line presents the plot for this choice of the bin size, the blue line is for roughly
ten times larger division $\Delta D=0.005$  while  green plot  presents  the histogram of prime pairs with $D$  divided into
bins of the size $\Delta D=10^{-1}$. These plots can be normalized by dividing all values by the maximal value
present in the histogram for a given bin size.

\begin{figure*}[t]
\begin{minipage}[t]{0.45\textwidth}
\includegraphics[width=\textwidth, angle=0]{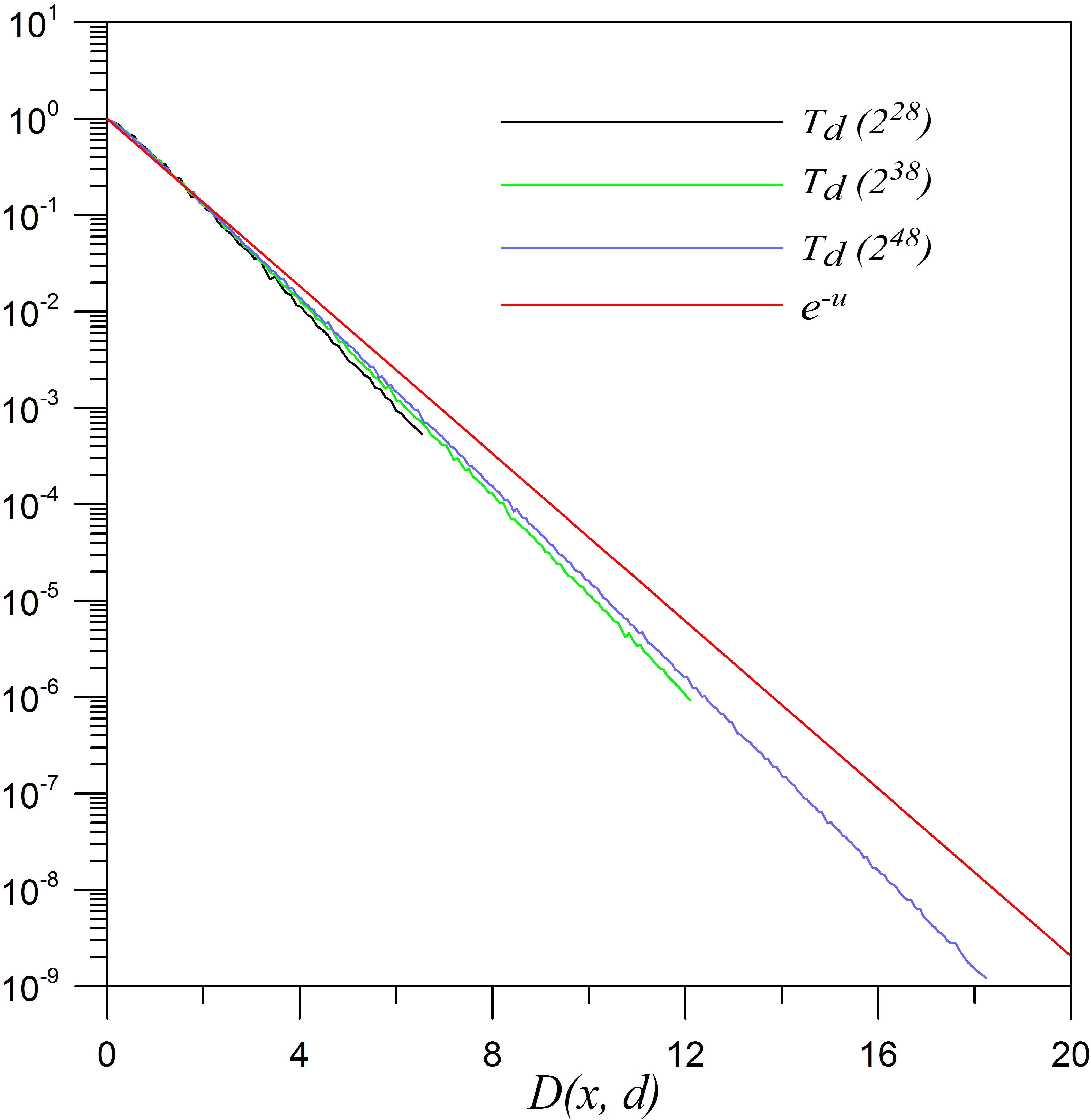}
\caption{Plots of $(\mathcal{D}(x, d), \mathcal{T}_d(x)),~(d=2, 4, \ldots)$ for $x=2^{28},2^{38},2^{48}$ and in red
the plot of $e^{-u}$. Only the points with $\tau_d(x)>1000$ were plotted to avoid
fluctuations at large $D(x, d)$ due to small values of $\tau_d(x)$ for large $d$.}
\label{fig-scaling}
\end{minipage}
\hspace{0.03\textwidth}
\begin{minipage}[t]{0.45\textwidth}
\centering
\includegraphics[width=\textwidth,  angle=0]{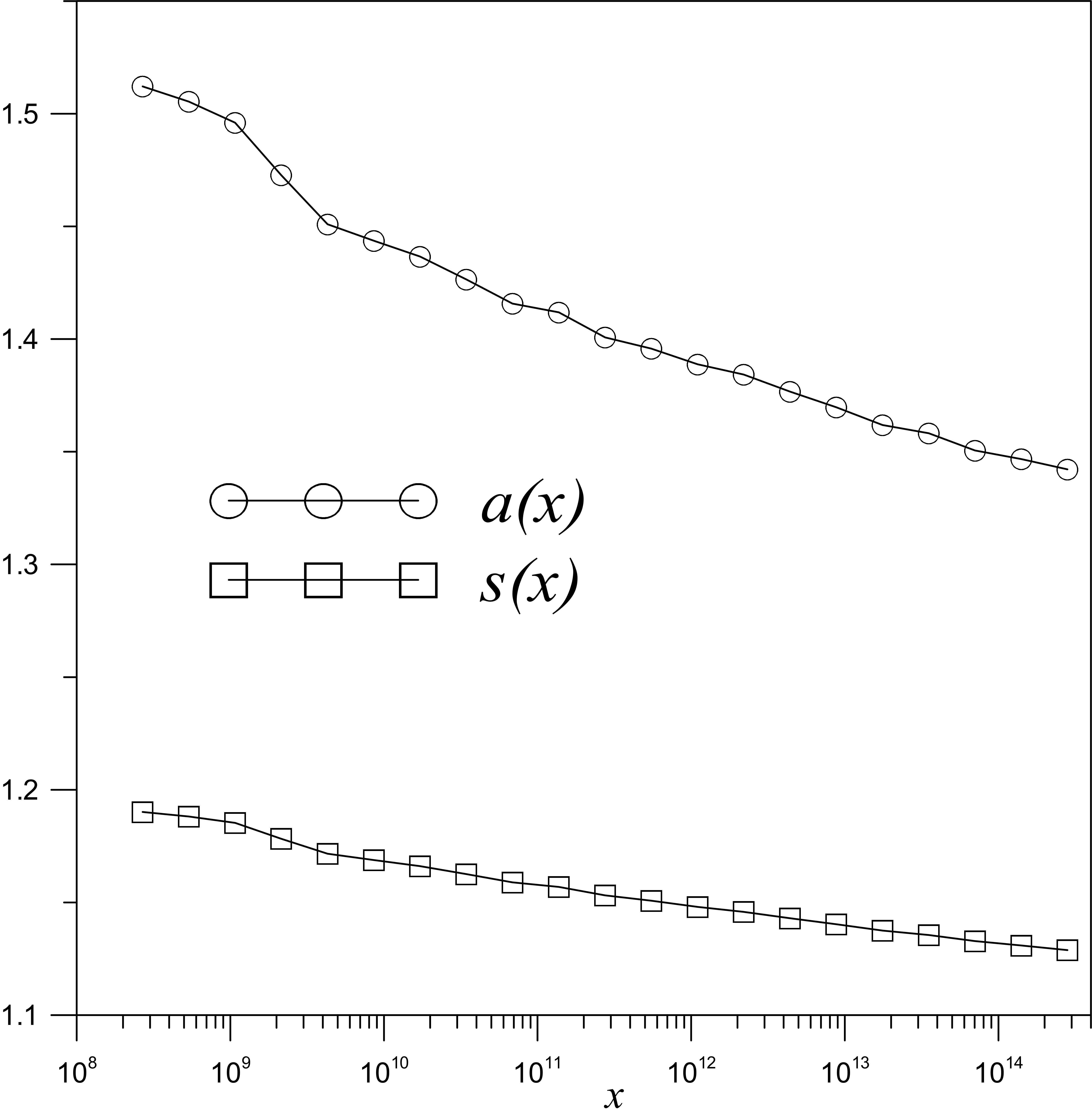}
\caption{ Plot of slopes $s(x)$ and prefactors $a(x)$ in the  dependence  $a(x)e^{-s(x)}$ obtained from fitting it to %   $(D(d, x), \ln(T_d(x)))$
$(\mathcal{D}(d, x), \ln(\mathcal{T}_d(x)))$ for $x=2^{28},2^{29}, \ldots, 2^{48}$.}
\label{fig-slopes}
\end{minipage}
\hspace{0.03\textwidth}
\vspace{0.6cm}
\end{figure*}

The explicit form of the equation (\ref{main}) allows us to define the unfolding in the following way:
Let us define the rescaled quantities:
\bee
\mathcal{T}_d(x)=\frac{x\tau_d(x)}{C_2 P(d) \pi^2(x)}, ~~~~~~\mathcal{D}(x, d)= \frac{d\pi(x)}{x}.
\label{eq-scaling}
\eee
The product $P(d)$ in the denominator of the first formula removes the oscillations and gives the analog of the histogram
free of size bin ambiguity.   The second  equation defines the proper unfolding for prime numbers.
Because $x/\pi(x)\approx \ln(x)$ is the mean distance between two consecutive primes  $\overline{d}\approx \ln(x)$     %  $<\!\! d \!\!>\approx \ln(x)$
up to $x$, we see that $\mathcal{D}(x, d)$ corresponds to the distances between  ``unfolded''  primes ---
normalized spacing between two consecutive primes is $\mathcal{D}(x, d)\approx d/\ln(x)$ and hence the mean value of
$\mathcal{D}(x,d)$ is simply 1. For large $x$  the quantity $\mathcal{D}(x, d)$ agrees with  expression (\ref{D_n-formula}) for large $p_n$:
$\mathcal{D}(p_n, d_n)\approx d_n/\ln(p_n)=D_n$ and hence values of $\mathcal{D}(x, d) \in [2/\ln(x),\, \ln(x)]$. From
the conjecture (\ref{main}) we expect that for each $x$  the points
$(\mathcal{D}(x, d), \mathcal{T}_d(x)), ~~d=2, 4, \ldots, G(x) $ should coincide --- the function
$\tau_d(x)$ displays scaling in the physical terminology. In Fig. \ref{fig-scaling}
we have plotted the points $(\mathcal{D}(x, d), \mathcal{T}_d(x))$ for $x=2^{28}, 2^{38}, 2^{48}$.
and indeed we affirm the tendency of all these curves to collapse into the universal one.
To make this plot we have used {\it exact} values of $\pi(x)$,  not any of the approximate formulas like $\Li(x)$:
from the definition of $\tau_d(x)$ it follows that $\pi(x)=\sum_d \tau_d(x)+1$ and it allowed us to calculate
from  $\tau_d(x)$ precise values of $\pi(x)$ for $x=2^{28}, 2^{38}, 2^{48}$.
If we denote $u=\mathcal{D}(x, d)$ then all these scaled functions should exhibit the pure exponential
decrease $e^{-u}$: Poisson distribution shown in red in Fig. \ref{fig-scaling}.
We have determined by the least square method  slope $s(x)$ and  prefactor  $a(x)$ of the fits
$a(x)e^{-s(x)u}$ to the linear parts of plots of $(\mathcal{D}(x, d), \ln(\mathcal{T}_d(x)))$. The results are
presented in Fig. \ref{fig-slopes}. The slope very slowly tend to 1: for over 6
orders of $x$ $s(x)$  changes  from $1.187$ to $1.136$  while the prefactor $a(x)$ drops from 1.512... to 1.273... .

Finally let us remark that there is no repulsion of small gaps between primes:  usually for GOE or GUE there is a prohibition
of small gaps between energy levels (in fact the number of gaps with $s=0$ is equal to zero), but for our case the smallest gap corresponds to
twins and it is believed that there is infinity of them. From (\ref{main}) it follows that the number of twins and cousins is roughly
a half of the number of primes separated by $d=6$.   In fact for all plots of $\tau_d(x)$ in Fig. \ref{fig-hist} $d=6$ is  the
highest point --- i.e.  it is most often occurring gap.  However in Fig. \ref{fig-hist}  local spikes appear at multiplicities of
$30=2\cdot 3\cdot 5$ and  at $d=210=2\cdot 3\cdot 5\cdot 7$,  where the product $P(d)$ has local maxima.  As $x$ increases the
slopes of plots of $\tau_d(x)$  decrease and at some value around $x\approx 10^{36}$ the peak at $d=30$ will be greater than that at
$d=6$. At much larger $x \approx 10^{428}$ the spike at $210$ will take  over $d=30$. It leads to the so called problem of
champions, i.e. most occurring  gap between consecutive primes, see  \cite{champions}. Thus primes are repelled in a very special
way: the most often occurring gaps are products of consecutive primes, but they become the ``champions'' at extremely
large values of $x$. For the unfolded according to eq.  (\ref{D_n-formula}) gaps $D_n$ (or eq.(\ref{eq-scaling}) and
quantities $\mathcal{D}$  as well)   there is no repelling:
the most common value of  $D_n$  is $2\pi(x)/x \approx 2/\ln(x)$ and it tends to zero with increasing $x$ ---behavior
typical for the Poisson  distribution.

Similar unfolding procedure has been used in dynamical systems e.g. in the stadium billiard were the
existence of strong oscillations due to bouncing ball orbits strongly
influence the spectral statistics $\Delta_3$ and to get a good agreement with the
Gaussian Orthogonal Ensemble (GOE)  predictions one has to perform unfolding which includes explicitly
the contribution of the bouncing ball periodic orbits (see \cite{Sieber-1993}).

%The explicit dependence on $x$ of the quantities $\tau_d(x)$  and ``unfolded level spacings'' $D(x, d)\sim d/\ln(x)$ shows
%that gaps between  primes are not described  by the stationary  Poisson distribution.
It is a common belief that the Poisson NNSD of the quantum energy levels is linked with
integrable  systems  with more than one degree of freedom.  In \cite{Crehan-1995} P. Crehan  has shown that for any
sequence of energy levels   obeying a certain growth law  ($|E_n|<e^{an+b}$, for some $a\in \mathbb{R}^+$, $b\in \mathbb{R}$),
there are infinitely many {\it classically  integrable}   Hamiltonians for
which the corresponding quantum spectrum coincides with this sequence. Because from PNT it follows, that the $n-$th prime
$p_n$ grows like $p_n\sim n \ln(n)$ the results of Crehan's paper can be applied and there exist classically integrable
hamiltonians whose spectrum coincides with prime numbers, see also \cite{Rosu-2003}.

\bigskip

\section{Spectral rigidity of prime numbers}
%\bigskip

In \cite{Mehta-Dyson-IV} several statistical measures to describe fluctuations in the energy levels $\{E_n\}$ of complex systems
were introduced. One which attracted much attention is the  spectral rigidity  $\Delta_3$. The  spectral rigidity
for arbitrary system with spectral staircase $N(E)$ is defined
as the averaged mean square deviation of the best local fit straight line $a\epsilon +b$ to the $N(E)$ on the interval $(x, x+L)$:
\bee
\Delta_3(x;L)=\frac{1}{L}\left \langle\min_{a, b} \int_0^L \left(N(x+\epsilon)-a\epsilon-b\right)^2 d\epsilon \right \rangle
\label{def-Delta_3}
\eee
The averaging  procedure  $\langle \cdot \rangle$ depends on the specific problem, e.g. for random matrices it is the mean value
from an ensemble of generated matrices or average over a set of atomic nuclei in real experiments, see e.g. \cite{Haq-et-al};
sometimes average over the initial point $x$ is applied.
There are in general  two ways of performing the operation $\min_{a, b}$, see the discussion in \cite{Mehta-Dyson-IV}.
One can calculate partial derivatives of r.h.s. of (\ref{def-Delta_3})
with respect to  $a$ and $b$, equate them to zero,  solve for $a, b$  and substitute solutions back to r.h.s. what
leads to the double integrals, see e.g.  \cite[Appendix II] {Pandey1979}. We will present here the procedure
for calculating  $\Delta_3$  in this way  devised by O. Bohigas and M.-J. Giannoni in \cite{Bohigas-Giannoni-1975} and \cite{Bohigas-1984}.
First the energies are unfolded  $E_N \rightarrow e_n$ using the smooth part $\overline{N}(E_n)$ of the staircase function,
see eq. (\ref{unfolding}). If the sequence of unfolded levels $e_1, e_2, \ldots, e_n$  falls  in the interval $(x, x+L)$ the
following explicit formula for $\Delta_3(x;L)$ is obtained:
\begin{widetext}
\bee
\Delta_3(x;L)=\frac{n^2}{16}-\frac{1}{L}\left(\sum_{k=1}^n \tilde{e_k} \right)^2 + \frac{3n}{2L^2}\sum_{k=1}^n \tilde{e_k}^2
-\frac{3}{L^4}\left(\sum_{k=1}^n \tilde{e_k}^2 \right)^2 + \frac{1}{L}\sum_{k=1}^n (n-2k+1)\tilde{e_k},
\label{Bohigas-Giannoni}
\eee
\end{widetext}
where  $\tilde{e_k}=e_k-(x+L/2)$. In the second approach  the  parameters
$a$ and $b$ are obtained by fitting the straight line $ax+b$ to the set of points $(x_1, y_1), (x_2, y_2), \ldots, (x_n, y_n)$
by the  least square method, i.e. the  partial derivatives of $\sum_{k=1}^n (y_k-ax_k-b)^2$ with respect to  $a$ and $b$
are calculated and put  equal to zero, what gives the very well known expressions:
\[
a=\frac{n\sum_{k=1}^n x_k y_k -\sum_{k=1}^n x_k \sum_{k=1}^n y_k}{n  \sum_{k=1}^n x_k^2 - (\sum_{k=1}^n x_k)^2}
\]

\[
b=\frac{1}{n} \sum_{k=1}^n \left( y_k - a x_k \right)
\]

\noindent In the case  of $\Delta_3(x;L)$ we have $x_k=E_k,  y_k=N(E_k)$.
The spectral rigidity obtained in this second way we will distinguish from (\ref{Bohigas-Giannoni}) by apostrophe
$\Delta_3'(x;L)$.  The formula for $\Delta_3'(x;L)$ in this approach and  adjusted
for our problem will be given below, see (\ref{moja_Delta}).

\begin{figure*}[t]
\begin{minipage}[t]{0.45\textwidth}
\centering
\includegraphics[width=\textwidth]{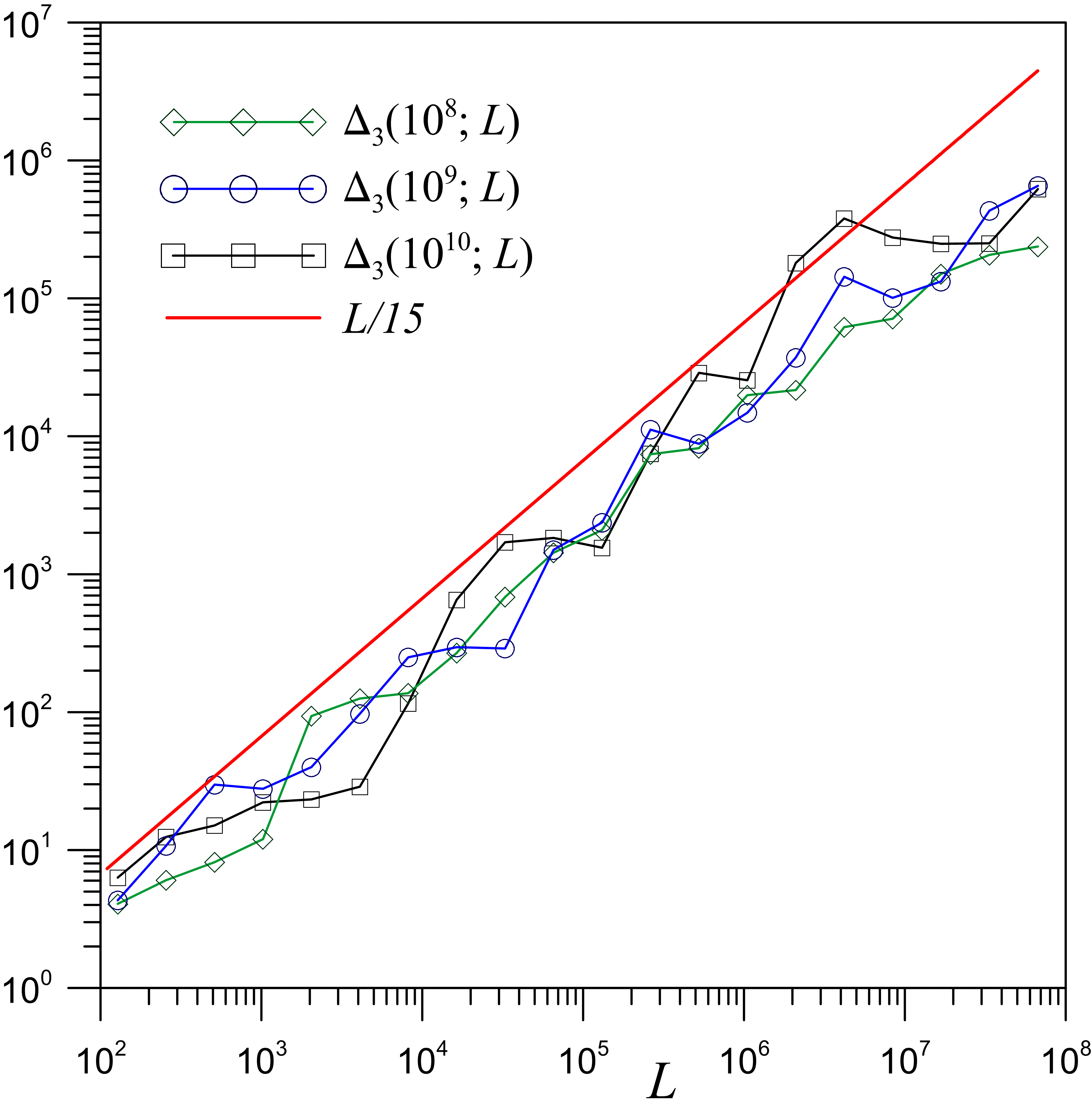}
\caption{ Plots of $\Delta_3(x; L)$ obtained from (\ref{Bohigas-Giannoni}) for $x=10^8, x=10^9$ and $x=10^{10}$
and $L=2^{7}=128, \ldots, 2^{26}=67108864$. }
\label{fig-Bohigas-Giannoni}
\end{minipage}
\hspace{0.03\textwidth}
\begin{minipage}[t]{0.45\textwidth}
\centering
\includegraphics[width=\textwidth, angle=0]{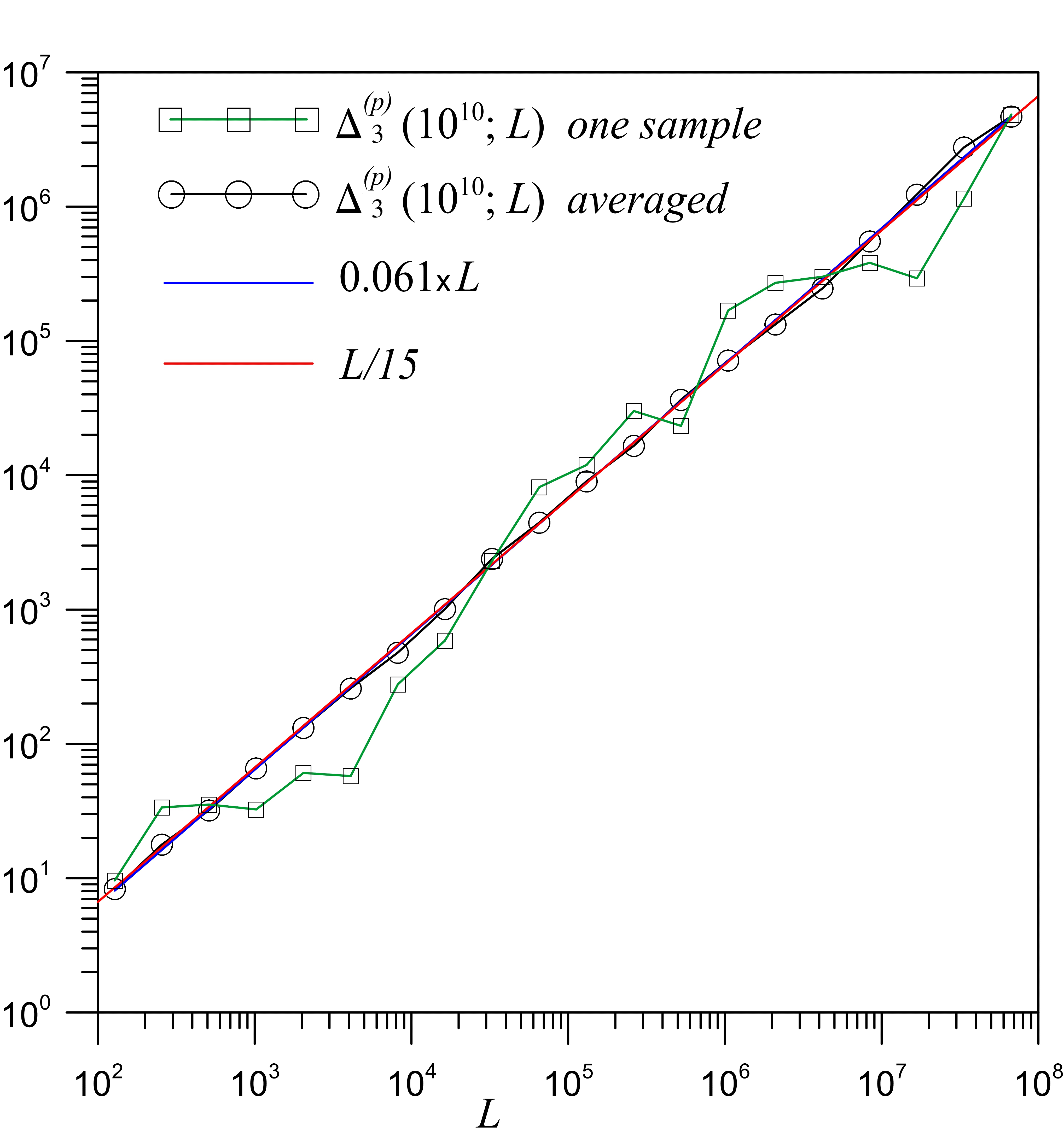}
\caption{The plot of $\Delta_3^{(p)}(L)$ for probabilistic primes for one particular realization of the ``artificial  primes''
in green (boxes) and averaged over 100 samples in black (circles). The last plot perfectly coincides with the red line
representing  $L/15$.  In blue is the fit $0.061L$ to circles plotted.}
\label{artificial}
\end{minipage}
\end{figure*}

Spectral rigidity for primes we  define  by (\ref{def-Delta_3}) with $\pi(x)$ instead of $N(x)$.
To use the formula (\ref{Bohigas-Giannoni}) the exact values of all primes are needed and we have used primes $p_n$
sufficient for calculation of $\Delta_3(x;L)$ for $x=10^8, 10^9$ and $10^{10}$.  To perform the unfolding $p_n \rightarrow r_n$
one can use in principle any  analytical formula  giving the number $\pi(x)$  of primes smaller than $x$, e.g. the one due
to Gauss $\pi(x)\sim x/\ln(x)$ or another one  given by the Prime Number Theorem (\ref{PNT}):  $r_n={\Li}(p_n)$.
The choice $ x/\ln(x)$  is  not a good one because  $\pi(x)-x/\ln(x)$ never changes the sign (see e.g. \cite[eq. 3.5]{Schoenfeld1962})
so there are no oscillations of  this  difference. Although  J.E. Littlewood has proved  in 1914 \cite{Littlewood},
that $\pi(x)-{\Li}(x)$  infinitely often changes the sign,  the lowest present day known estimate for the first sign
change of $\pi(x)-{\Li}(x)$  is around $10^{316}$, see \cite{Bays} and \cite{Demichel2}, hence  in the available for computers
range there are  no fluctuation of $\pi(x)-{\Li}(x)$ around zero but a steady growth of the function ${\Li}(x)-\pi(x)$.
In the famous paper \cite{Riemann}  B. Riemann has given the {\it exact} formula for $\pi(x)$:
\begin{widetext}
\bee
\pi(x)=\sum_{k=1}^\infty \frac{\mu(k)}{k}\left({\Li}(x^{\frac{1}{k}})-\sum_{\rho} {\Li}(x^\frac{\rho}{k})  +\int_{x^{1/k}}^{\infty } \frac{1}{u \left(u^2-1\right)\ln(u)} \, du\right)
\label{Riemann}
\eee
\end{widetext}
\noindent where  $\mu(n)$ is the M{\"o}bius function:
\[
\mu(n) \,=\,
\left\{
\begin{array}{ll}
1 & \mbox {for  $ n =1 $} \\
0 & \mbox {when  $p^2|n$}\\
(-1)^r & \mbox{\rm when}~ n=p_1 p_2 \ldots p_r.
\end{array}
\right.
\]
The sum over $\rho$ runs over nontrivial zeros  of the Riemann $\zeta(s)$ function $\zeta(\rho)=0$ and the last integral contains
contribution  from trivial zeros $-2m$ of  zeta: $\zeta(-2m)=0, ~m=1,2,3, \ldots$. If the Riemann Hypothesis is true then  for
all nontrivial zeros $\Re(\rho)=\frac{1}{2}$  and the contribution to  the sum  over $k$ in (\ref{Riemann}) is dominated by
the first term,  what leads to the  following approximation  to $\pi(x)$:
\bee
R(x)= \sum_{k=1}^\infty \frac{\mu(k)}{k}{\Li}(x^{\frac{1}{k}}).
\label{def-R}
\eee
The difference $\pi(x)- R(x) $ changes the sign already at $x$ as low as
$x\in [2, 100]$,  see e.g. tables  obtained by T. R. Nicely in \cite{web-primes} and up to $10^{14}$ there are over
50 millions of sign change of $\pi(x)-R(x)$ \cite{Kotnik2013},  however on  average the behavior of both differences
$\pi(x)-{\Li}(x)$   and $\pi(x)-R(x)$ seems to be the same \cite{Kotnik2008}. The above function $R(x)$ can be obtained,
without the need of calculating  the logarithmical integral ${\Li}(x)$,  from the series obtained  by J.P.  Gram,
see e.g.  \cite[p.51] {Riesel}:
\bee
R(x)=1+\sum_{m=1}^\infty \frac{\ln^m(x)}{m m! \zeta(m+1)}
\label{eq-Gram}
\eee
Hence we have made the unfolding of primes according to the rule
\bee
r_n=R(p_n).
\label{unfolding-R}
\eee
At this point let us remark that from (\ref{Li-series})
and (\ref{eq-Gram}) we see that because $\zeta(m)\to 1$ for $m \to \infty$ very quickly (e.g. $\zeta(4)=\pi^4/90=1.082323\dots,
~\zeta(6)=\pi^6/945=1.017343\ldots $) for large $x$  the functions $\Li(x)$ and  $R(x)$
should differ  by roughly $\ln \ln(x)$ and this quantity can be discarded in comparison with values of series involving powers
of $\ln(x)$ present in (\ref{Li-series}) and (\ref{eq-Gram}).
Indeed,   from (\ref{def-R}) it follows using the first term from asymptotic expansion (\ref{PNT})
that for large $x$ the approximate relation $R(x)/\Li(x)=1-1/\sqrt{x}$ holds.
Thus for large $x$ the particular form of unfolding (\ref{unfolding-PNT}) or
(\ref{unfolding-R})  should be irrelevant, despite the fact that $\pi(x)-\Li(x)$ changes the sign first time somewhere in the vicinity
of $x=10^{316}$ while $\pi(x)-R(x)$ changes the sign already for $x$ between $10$ and $20$, see tables of Nicely \cite{web-primes}.

We will present the plots of $\Delta_3(x;L)$ for three values of $x$: $10^8, 10^9$ and $10^{10}$.  The values of primes for which
the unfolded variables begin to fall into the intervals $(10^8, 10^8+L), (10^9, 10^9+L), (10^{10}, 10^{10}+L),$ are accordingly
$ ~2038076627, ~22801797631, ~252097715777$: $R(2038076627)=10^8+1.8496\ldots, R(22801797631)=10^9+2.3178\ldots, R(252097715777)= 10^{10}+0.0024\ldots$.
As there seems to be no clear relation between
the  values of  $L$ in comparison with chosen $x$ we have used the wide range of values of $L$:  we have  calculated from
(\ref{Bohigas-Giannoni}) spectral rigidity for values $L=2^7=128,   \ldots, 2^{26}=67108864$.
The results are presented in Fig. \ref{fig-Bohigas-Giannoni}. It is well known that for stationary Poisson ensemble
$\Delta_3(x;L)=\frac{L}{15}$,  see e.g. \cite[eq.(61)]{Mehta-Dyson-IV} or \cite[Appendix II]{Pandey1979},
and on the Fig. \ref{fig-Bohigas-Giannoni} this theoretical prediction is plotted in blue.
The obtained plots of $\Delta_3(x;L)$ seem to tend to the line $L/15$ with increasing $x$.
For primes  there is no natural averaging procedure present in the definition (\ref{def-Delta_3})
and in Fig. \ref{fig-Bohigas-Giannoni}  prominent fluctuations are seen.  To simulate the averaging we have performed the
following ``Monte Carlo'' experiment for $x=10^{10}$. From the  PNT in the form $\pi(k)\sim k/\ln(k)$ it follows that the
chance that randomly chosen large integer $k$ should be a  prime is  $1/\ln(k)$.  Such a probabilistic model for primes
was created by H. Cramer in  the 1930's  \cite{Cramer}.  We have started to test if a given natural $k$ number is  the probabilistic
``artificial '' prime  from the first $k_0$ for which  $R(k_0)>10^{10}$, i.e. for $k_0=252097715777$ for which $R(k_0)=10^{10}+0.00241\ldots$.
The natural number $k>k_0$ (even the even numbers were allowed --- when even numbers are skipped the probability of odd number $k$
to be a ``prime'' should be $2/\ln(k)$) was accepted to be a  ``probabilistic'' prime if $1/\ln(k)$ was larger than the
uniformly generated  from the interval $(0,1)$  random number \verb"random": \verb"random"$\,\, <\! 1/\ln(k)$. For such a ``prime'' $k$
the unfolding was performed using the equation $r'_k=R(k)$.  The random  drawing of ``primes''  was  continued  until the
unfolded ``prime'' was larger than $x+L$ for $L=128, \ldots, 2^{26}$.  For the set of such generated unfolded quantities
in the intervals $(x, x+L)$  the ``artificial''  spectral rigidity  $\Delta_3^{(p)}(x;L)$  was  calculated using
(\ref{Bohigas-Giannoni}).  The result of this procedure is plotted  in green in Fig. \ref{artificial} and there are
fluctuations seen resembling those present  in Fig. \ref{fig-Bohigas-Giannoni} for ``true'' primes.
But now we can generate many independent sets of the  artificial probabilistic primes.  We have repeated this
procedure  100 times and the averaged over all these  samples spectral rigidity $\Delta_3^{(p)}(x; L)$ is presented in Fig. \ref{artificial}
in black. Now the fluctuations disappeared  and the obtained plot follows perfectly the predicted dependence $L/15$.
This allows us to claim that  the spectral rigidity for prime numbers unfolded via the Riemann function $R(x)$  is the same
as for Poisson  statistics  (we have checked that the same result is obtained for unfolding with $\Li(x)$ as in eq. (\ref{unfolding-PNT})).
Let us mention that usually saturation of  $\Delta_3$  is observed in physical systems,
i.e. after the initial dependence resembling $L/15$  spectral rigidity  stops to increase and is constant for large $L$, see
e.g.   \cite{Berry-1985} or \cite{Casati-et-al-1985}. However our system is infinite  and there is no departure from
straight line $L/15$.

Next  we will present spectral rigidity for  second method of minimizing the r.h.s of (\ref{def-Delta_3}) over $a, b$, namely
determination of $a, b$ by the least  square method.  In the case of primes numbers, for large $x$,  the smooth part
of staircase $\pi(x)$ given by $x/\ln(x)$ is almost linear in the interval $(x, x+L)$ as the denominator  changes from
$\ln(x)$  to  $\ln(x+L)=\ln(x)+L/x+\ldots$  what  for $x\gg L \gg 1$ again is $\ln(x)$.
There are a few websites \cite{web-primes}   offering the tables of values of $\pi(x)$ (as well as other number theoretic
functions). In these  data files the values of $\pi(x)$ are tabulated with different step size of $x$, the best resolution
is at the  A. V. Kulsha's page: the file pi.txt of the size  421MB  contains counts of $\pi(x)$ with a step of $10^9$  from
$x=10^9$  to $x=2.5\times 10^{16}$. Now we will give the formula for calculating the integral appearing in the definition
of $\Delta_3'(x;L)$:
\bee
\mathcal{I}(x; L)=\int_0^L \left(\pi(x+\epsilon)-a\epsilon-b\right)^2 d\epsilon  %   \int_0^L \left(\pi(x+u)-au-b\right)^2 du
\label{calka}
\eee
appropriate for our data.  Let us assume, that the values of $\pi(x)$ in the integral (\ref{calka}) are known with the resolution $h$:  $y_k=\pi(x+kh)$;
hence we  assume  that $\pi(x)$ is constant on the  intervals $(kh,  (k+1)h)$ (in fact $\pi(x)$ is constant only between
two consecutive  primes). We regard this sampling of $\pi(x)$ with different steps $h$ as the averaging procedure hidden in the
angle bracket in (\ref{def-Delta_3}) ---  taking values of $\pi(x)$ at all consecutive  primes would  introduce
fluctuations. The  combination $\pi(x+\epsilon)-a\epsilon-b$
is the linear function on the intervals  $(kh,  (k+1)h)$ and we can write (we assume here that $L$ is the integer multiple of $h$):
\begin{eqnarray*}
\mathcal{I}(x; L)=\int_0^L \left(\pi(x+\epsilon)-a\epsilon-b\right)^2 d\epsilon = \\
 \sum_{k=0}^{L/h - 1} \int_{kh}^{(k+1)h} \left(y_k-a\epsilon-b\right)^2 d\epsilon.
\end{eqnarray*}

Performing elementary integration we obtain:
\begin{widetext}
\bee
\Delta_3'(x; L)=b^2+abL+\frac{a^2L^2}{3} + \frac{1}{L}\sum_{k=0}^{L/h-1} y_k(y_k-2b)h-ay_k(2k+1)h^2
\label{moja_Delta}
\eee
\end{widetext}
It should be noted, that parameters  $a$ and $b$ in eq. (\ref{moja_Delta}) obtained from  fitting  $a\epsilon + b$  to  points
$\pi(x+\epsilon),  0\leq \epsilon \leq L$,  by least-square method are functions of $L$ and $x$, see below (\ref{a-b}).

\begin{figure}[h]
\includegraphics[width=0.3\textwidth]{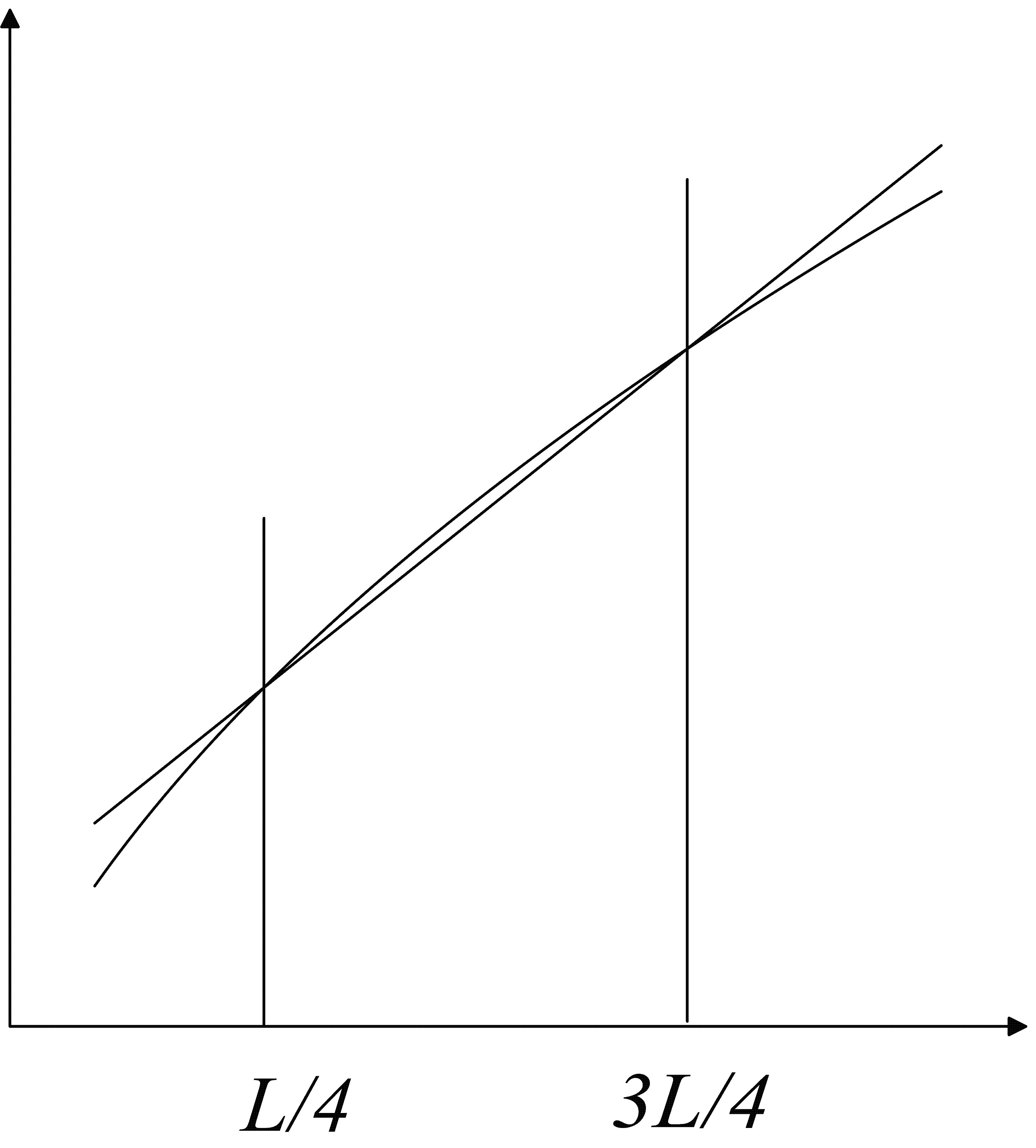}
\caption{ The illustration of the experimental fact that the straight line best
fitting  $(x+\epsilon)/\ln(x+\epsilon)$ on the interval $\epsilon \in (0, L)$  crosses it at $\epsilon=L/4$ and $\epsilon=3L/4$.}
\label{fig-LSQM}
\end{figure}

\begin{figure*}[t]
\begin{minipage}[t]{0.45\textwidth}
\centering
\includegraphics[width=\textwidth]{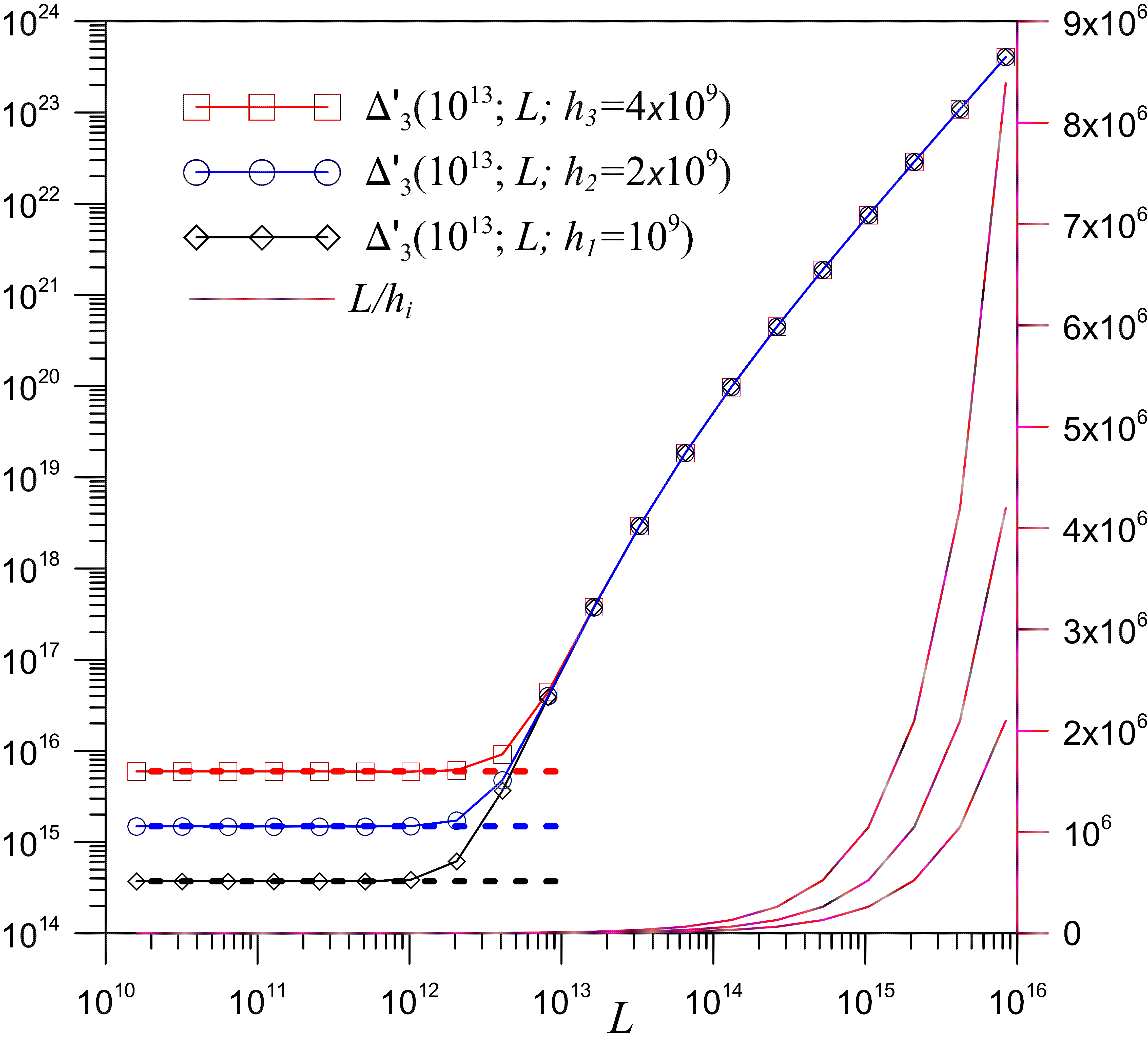}
\caption{Plots of $\Delta_3'(x_1; L; h)$  for $x_1=10^{13}$ and  $L=16h_1=1.6\times 10^{10}, 32h_1=3.2\times 10^{10},
\ldots 2^{23}h_1=8.388608\times 10^{15}$ and three values of $h_1=10^9$,  (black), $h_2=2h_1$  (blue)  and $ h_3=4h_1$ (red).
On the right in regal red are plotted values of the number of terms $L/h_i-1$
in the  sum (\ref{moja_Delta}) and the right $y$  axis also in  regal red is for these numbers. The dashed lines represent
values of (\ref{delta-od-h}).   The coincidence of $\Delta_3'(x_1; L)$'s for all $h_i$ starts at approximately
$L=2^{14}h_1= 1.6384\times 10^{13}$.}
\label{fig-delta1}
\end{minipage}
\hspace{0.03\textwidth}
\begin{minipage}[t]{0.45\textwidth}
\centering
\includegraphics[width=\textwidth]{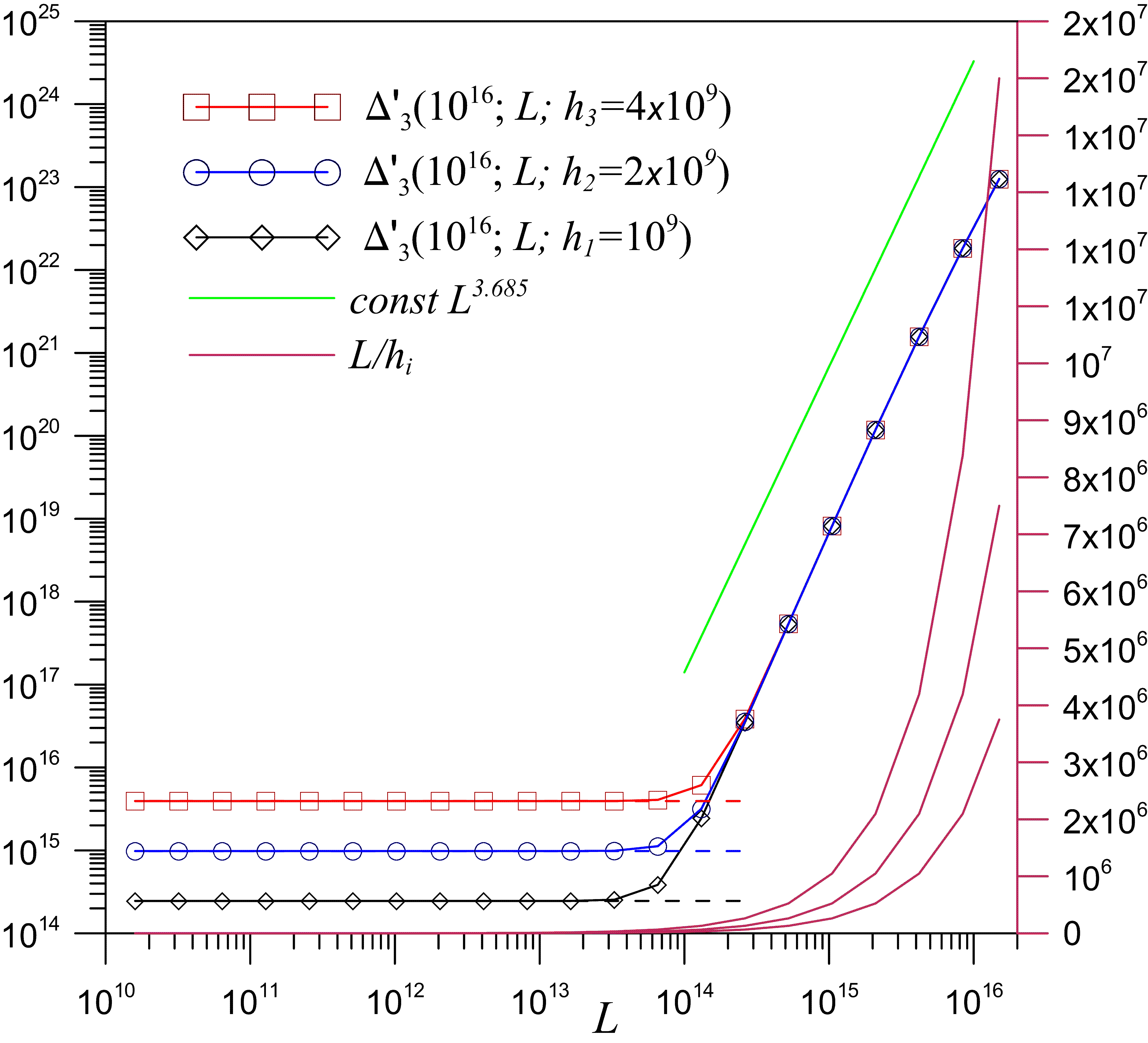}
\caption{Plots of $\Delta_3'(x_2; L; h)$  for  $x_1=10^{16}$ and  $L=16h_1=1.6\times 10^{10}, 32h_1=3.2\times 10^{10},
\ldots 2^{23}h_1=8.388608\times 10^{15}$ and additionally for  $L=1.5\times 10^{16}$  and three values of  $h_1=10^9$,
(black), $h_2=2h_1$  (blue)  and $ h_3=4h_1$ (red). On the right in regal red are plotted values of the number of terms $L/h_i-1$
in the  sum (\ref{moja_Delta}) and the right $y$ axis  also in  regal red is for these numbers. The dashed lines represent
values of (\ref{delta-od-h}).  The coincidence of $\Delta_3'(x_2; L)$'s for all $h_i$ starts at approximately
$L=2^{18}h_1= 2.62144\times 10^{14}$  and follows practically  power-like increase given by equation
$ 2.5695\times 10^{-37}L^{3.685}$ --- the  green line presents this equation multiplied by 100, however we expect bending of
$\Delta_3'(x_2; L; h)$ for $L>x_2$ similar to the one seen in Fig. \ref{fig-delta1}. }
\label{fig-delta2}
\end{minipage}
\end{figure*}

The value of  $\Delta_3'(x;L)$ given by (\ref{moja_Delta})
should not depend on $h$.  To test this presumption we have calculated $\Delta_3'(x;L)$ for $x_1=10^{13}$ and $x_2=10^{16}$
and for $h_1=10^9, h_2=2\times 10^9, h_3=4\times 10^9$. We have chosen  the  following sequence
of values of the length of intervals $L=16h_1=1.6\times 10^{10}, 32h_1=3.2\times 10^{10}, \ldots 2^{23}h_1=8.388608\times 10^{15}$  for both $x_1, x_2$ and
additionally   $L=1.5\times 10^{16}$ for $x_2=10^{16}$. It means that the number of terms in the sum in
(\ref{moja_Delta}) was $2^3, 2^4, \ldots, 2^{22}=4194304$ for $h_2$ and   $2^2, 2^4, \ldots, 2^{21}=2097152$  for $h_3$
appropriately.   For each $L$ the parameters $a$ and $b$  were fitted by the least-square method to the points
$(x_k=x+kh, ~y_k=\pi(x+kh)), ~k=0, 1,  \ldots, L/h-1$.
In Figures \ref{fig-delta1} and \ref{fig-delta2} we present the results.  Two types of behaviors are  seen
in these figures:  the constant in  $L$ values of $\Delta_3'$ depending on $h$ and the collapse of plots of  $\Delta_3'$
for all $h$ when  the  increase of $\Delta_3'$ with $L$ begins. It seems that to get rid of dependence on  $h$ the sufficiently
large number $L/h$ of terms  in  the  formula (\ref{moja_Delta}) has to be summed up.  The inspection of data shows, that to have
the independence of $\Delta_3'$ on $h$ a few thousands of terms in the sum in (\ref{moja_Delta}) are sufficient (for largest $L$
there are millions of terms in this sum, see plots  in royal red in Fig. \ref{fig-delta1} and \ref{fig-delta2}).
It is possible to find heuristically the values of the constant in $L$ parts of $\Delta_3'$.  To find the analytical expressions
for $a$ and $b$  we consider the smooth part
of $\pi(x)$ given by $(x+\epsilon)/\ln(x+\epsilon)$ and the straight line $a\epsilon+b$ obtained by best fitting to the values of
$(x+kh)/\ln(x+kh)$.
The experiments show, that the fits cross $(x+\epsilon)/\ln(x+\epsilon)$ on the interval $\epsilon\in (0, L)$ roughly at
$\epsilon=L/4$ and $\epsilon=3L/4$,  see  Fig. \ref{fig-LSQM},  thus from  $(x+L/4)/\ln(x+L/4)=aL/4+b$
and $(x+3L/4)/\ln(x+3L/4)=a3L/4+b$  we get
\begin{widetext}
\begin{eqnarray}
a=\frac{2}{L} \left(\frac{x+3L/4}{\ln(x+3L/4)}-\frac{x+L/4}{\ln(x+L/4)}\right) =
\frac{1}{\ln(x)} -\frac{L}{2 x\ln^2(x)}+ { ~terms  ~~} \frac{ 1}{x^2}{ ~~or ~~ higher } \\%    \frac{13L^2}{96 x^2 \ln^2(x)} + \ldots\\
b=\frac{x+L/4}{\ln(x+L/4)} - aL/4=\frac{x}{\ln(x)}-\frac{L}{4 \ln^2(x)}+{~terms ~~ } \frac{ 1}{x} {~~ or ~~ higher } %  \frac{3L^2}{32 x \ln^2(x)} + \ldots
\label{a-b}
\end{eqnarray}
Using $y_k=(x+kh)/\ln(x+kh)\approx (x+kh)/\ln(x)-(kh)^2/2x\ln^2(x)$ we obtain in  (\ref{moja_Delta})  sums over $k$ which
can be calculated exactly and retaining the leading terms gives:
\bee
\Delta_3'(x;L)=\frac{h^2 }{3\ln^2(x)}-\frac{hL }{4\ln^3(x)}+\left({~ terms~~  } \frac{ 1}{x} { ~~or~~ higher }\right)
\label{delta-szereg}
\eee
\end{widetext}
Because $\Delta_3'(x; L)>0$ we have from above $h^2/ 3\ln^2(x)>hL/4\ln^3(x)$, i.e. $L<4h\ln(x)/3$, what for $x=10^{13}$ gives $L<40h$.
Surprisingly  the first term in (\ref{delta-szereg}), not depending  on  $L$ but being the function of $x$, gives the expression
\bee
\Delta_3'(x;L; h)=\frac{h^2 }{3\ln^2(x)}+\ldots
\label{delta-od-h}
\eee

which works very well even for $L=1024h$ for $x_1=10^{13}$ and $L=8192h$ for $x_2=10^{16}$, as it is seen in Figures
\ref{fig-delta1} and \ref{fig-delta2}, where the predicted values $h^2/ 3\ln^2(x)$  are plotted by dashed lines together
with the plots of $\Delta_3'(x;L; h)$ obtained  from (\ref{moja_Delta}).  In fact this agreement is astonishing:
e.g.  all $\Delta_3'(10^{16};L; h_1)$ for initial 11 values of $L$   have first three digits the same:
$2.455\ldots \times 10^{14}$  while (\ref{delta-od-h})  predicts $2.45588\ldots \times 10^{14}$.
In Fig. \ref{fig-delta1} we were able to make the plot for $L$ up to almost $10^3x_1$, while in Fig. \ref{fig-delta2}
the largest $L$ is smaller than $x_2$, thus we expect bending of $\Delta_3'(x_2;L; h)$ for larger $L$, similar to the
behavior of $\Delta_3'(x_1;L; h)$ on Fig. \ref{fig-delta1}.  In the plots of $\Delta_3'$ we see  the crossover
at value  $L^*$ above  which the steeper increase of spectral rigidities begins and  this  dependence is  $L^{\gamma}$,
with   $\gamma\approx 3.1$. Heuristically existence of this crossover
can be justified by the following reasoning:  for moderate values of $L$ the straight line $a\epsilon+b$  approximates
$\pi(x+\epsilon)$ quite well  leading to the small values of the integral $\int_x^{x+L} (\pi(x+\epsilon)-a\epsilon-b)^2
d\epsilon$,  while for larger $L$  the discrepancy between $\pi(x+\epsilon)$ and the straight line increases leading
to  larger values of $\Delta_3$. The spectral rigidity  calculated in second way displays different behavior than
$\Delta_3(x;L)$ obtained in the first manner.  Let us remark at this point that the proof of  $\Delta_3(x;L)=L/15$
for the Poisson  ensemble was obtained in \cite[Appendix II]{Pandey1979} only for the first method of minimalization
over $a$ and $b$ in  (\ref{def-Delta_3}).

\section{Conclusions}  In this paper we have treated prime  numbers as energy levels and we  applied the physical methods
used to study spectra of quantum systems to the description of distribution of prime numbers. We  presented large numerical
data (up to  $x=2.814 \ldots \times 10^{14}$) in support of the formula (\ref{main})  for NNSD between consecutive  primes.
It was also possible to obtain analytical formula (\ref{d_max}) for the maximal difference between two adjacent primes smaller
than $x$. The case of primes numbers gives the rare opportunity to calculate spectral  rigidity  $\Delta_3(x;L)$  for the wide
range of $x$ and $L$  --- for real physical systems usually only hundreds (nuclei), thousands or hundreds of thousands
(e.g. billiards)   energies are known. As the main result of this paper we regard the scaling relations (\ref{eq-scaling}) and
apparently  the first in the literature  attempt  to calculate spectral rigidity $\Delta_3(x; L)$ for prime numbers.
We have proposed the method to average the spectral rigidity over realizations of  probabilistic
primes and after sampling  over 100 sets of   ``artificial''  primes we have obtained  perfect $L/15$ dependence.
The obtained results confirm that the primes follow the Poisson distribution.  This averaging shows that the  spectral rigidity
does not depend on   peculiarities of the primes, but on the probability $1/\ln(k)$ of the number $k$ to be prime.
All the  above analysis can be repeated for subsets of prime numbers, for example for the twin primes (both $p$ and $p+2$ are prime),
cousin primes (both $p$ and $p+4$ are prime)  or the primes of the form $4k^2+1$; in the latter case the ``energy levels''
are the values of $k$ for which $4k^2+1$ is prime.

%\vspace{0.8cm}

\section{Acknowledgments}

I would like to thank  Prof. Marek Ku{\'s},  Prof. Jonathan Sondow  and Prof. Karol {\.Z}yczkowski for  comments and remarks.
I would like also to thank anonymous referees for useful comments and suggestions.

\end{document}